\documentclass[12pt]{amsart}
\usepackage{amscd,amssymb,amsmath,enumerate}
\usepackage{amsmath,amsthm}
\usepackage[all]{xy}
\usepackage{graphicx}
\usepackage{graphics}
\usepackage{latexsym}
\usepackage[all]{xy}
\usepackage{psfrag}
\xyoption{matrix} \xyoption{arrow}
\usepackage{parskip}
\usepackage{tikz}
\usepackage[justification=centering]{caption}

\newtheorem{proposition}{Proposition}
\newtheorem{theorem}[proposition]{Theorem}

\newtheorem{definition}[proposition]{{Definition}}
\newenvironment{defn}{\begin{definition} \rm}{\end{definition}}

\newtheorem{remark}[proposition]{{Remark}}

\newcommand{\G}{{ G}}

\begin{document}

\date{today}
\title[Skew-gentle algebras and orbifolds]
{Skew-gentle algebras and orbifolds}

\author[Schroll]{Sibylle Schroll}
\address{Sibylle Schroll\\
Department of Mathematics \\
University of Leicester \\
University Road  \\
Leicester LE1 7RH, UK
}
\email{schroll@le.ac.uk}



\setcounter{section}{1}
\keywords{}
\thanks{}

\begin{abstract}  This is an extended abstract of my talk at the Oberwolfach-Workshop "Representation Theory of Finite-Dimensional Algebras" (19-24 January 2020). It introduces  a geometric orbifold model for the bounded derived category of a skew-gentle algebra. This is joint work with Daniel Labardini-Fragoso and Yadira Valdivieso-D\'iaz. 
\end{abstract}
\date{\today}
\maketitle

Originating in cluster theory, geometric surface models have been instrumental in connecting representation theory with other areas of mathematics such as homological mirror mirror symmetry, see for example \cite{HKK, LP} and \cite{DJL}. 
In this extended abstract based on \cite{LSV},   we give an orbifold model for the bounded  derived category $D^b(A)$ of a skew-gentle algebra $A$ which encodes the indecomposable objects in $D^b(A)$
in terms of graded curves in the orbifold.

Skew-gentle algebras were introduced in \cite{GdP} in the context of the study of the Auslander Reiten theory of  clannish algebras. They are closely linked to the well-studied class of gentle algebras. 

A quiver is a  quadruple $Q= (Q_0, Q_1, s, t)$ consisting of a pair of finite sets, the vertex set $Q_0$ and 
the arrow set $Q_1$, and two maps $s, t: Q_1 \to Q_0$. We think of $Q$ as a directed graph and of 
elements in $Q_1$ as the arrows, that is, directed edges  $s(a) \to t(a)$, for $a \in Q_1$. 
Throughout let $K$ be an algebraically closed field.

\begin{defn} A $K$-algebra $A$ is skew-gentle if $A$ is Morita equivalent to an algebra $KQ/I$ where 
 $Q_1  = Q'_1 \cup S$ with $S \subset \{ a \in Q_1 \mid s(a) = t(a) \}$, $I = I' \cup \{\varepsilon^2 - \varepsilon \mid \varepsilon \in S\}$, and where
 $KQ' / I'$ is a locally gentle algebra with $Q' = (Q_0, Q'_1)$. Furthermore, if 
 $\varepsilon \in S$ then $s(\varepsilon)$, as a vertex in $Q'$, is a single arrow source, a single arrow sink or there exist exactly one arrow $a\in Q_1$  and one arrow $b \in Q_1$ with $s(b) = t(a)$ and   $ab \in I'$.  
 We call the elements of $S$ special loops.  
\end{defn}

We note that in the above definition even if the skew-gentle algebra $KQ/I$ is finite dimensional, the ideal $I$ is not necessarily admissible. We refer to \cite{GdP} for an admissible presentation $KQ^{sg}/I^{sg}$ isomorphic to $KQ/I$. 

Finite dimensional skew-gentle algebras are tame \cite{CB} and derived tame \cite{BMM, BD}. In this paper we will use the description of the indecomposable objects in $D^b(A)$, for $A$ finite dimensional skew-gentle, in terms of generalised homotopy strings and bands as given in \cite{BMM}. 

Given a  skew-gentle algebra $KQ/I$ with underlying gentle algebra $KQ'/I'$, we will now construct an orbifold dissection into generalised polygons. For this, recall from \cite{OPS} that gentle algebras (up to isomorphism) are in bijection with surface dissections into a special set of polygons (up to homeomorphism). Let $(\Sigma', \G',  \G'^*, M', P^*)$ be the surface dissection corresponding to $KQ'/I'$ as defined in \cite[\S1]{OPS}, where $\Sigma'$ is an oriented  surface with boundary and marked points, $\G'$ is a dissection of $\Sigma'$ into polygons, $\G'^*$ is a graph dual to $\G'$, $M'$ are those vertices of $\G'$ which are in the boundary of $\Sigma'$ and $P^*$ are those vertices of $\G'^*$ which are in the interior of $\Sigma'$.

We recall that the edges in $\G'$ correspond to the vertices in $Q'_0 = Q_0$. 
Suppose an edge $g' \in \G'$ corresponds to a vertex $v$ in $Q'_0$ which is the start of a special loop. Then $g'$ cuts out a digon in $(\Sigma', \G')$. We now describe a local replacement operation resembling in each step two consecutive Whitehead moves: namely, we contract the special edge $g'$ and expand it in the orthogonal direction producing in the process an edge g connecting a new marked point on the boundary and an orbifold point $\omega$ of order 2  as illustrated locally in the following  example: 

 \begin{figure}[ht!]
        \centering
        \begin{tikzpicture}
        \filldraw (0,0) circle (1pt) node[above] {\tiny$X$};
        \filldraw (2,0) circle (1pt) node[above] {\tiny$Y$};
        \filldraw (2.2,-1) circle (1pt);
        \filldraw (-0.4,-1) circle (1pt);
        \draw (0,0) ..controls (1,-0.5) ..(2,0);
        \draw[dashed] (-0.7,0.06) ..controls (1, -0.1).. (2.7,0.06);
        \draw (0,0) to[bend left] (-0.4,-1);
        \draw (-0.4,-1) to[bend left] (0.3,-1.7);
        \draw (2,0) to[bend right] (2.2, -1);
        \draw (2.2,-1) to[bend right] (1.5,-1.7);
        \node at (0.7,-1) {\tiny$$};
        \node at (1,-0.55) {\tiny{$g'$}};
        \draw (0,0) to[bend left] (-1,-0.3);
        \draw (2,0) to[bend right] (3,-0.3);
        \node at (-0.5,-0.5) {\tiny{$\dots$}};
        \node at (2.3,-0.5) {\tiny{$\dots$}};
        \node at (1,-1.7) {\tiny$\dots$};
        \node at (-1.5,0) {$\G' :$};
        \filldraw (7,0) circle (1pt) node[above] {\tiny$X=Y$};
        \filldraw (6,-0.5) circle (1pt);
        \filldraw (8,-0.5) circle (1pt);
        \draw[dashed] (5.5,0.45) to [bend right] (8.5,0.45);
        \draw (7,0) to[bend left] (5.5,0.2);
        \draw (7,0) to[bend right] (8.5,0.2);
        \draw (7,0) to[bend left] (6,-0.5);
        \draw (7,0) to[bend right] (8,-0.5);
        \draw (6,-0.5) to[bend left] (6.3, -1.3);
        \draw (8,-0.5) to[bend right] (7.7,-1.3);
        \node at (7,-1.5) {$\dots$};
        \draw (7,0) to (7,-1);
        \node at (7,-1) {$\times$};
        \node at (7.1, -1.2) {\tiny$\omega$};
        \node at (6.5,-0.7) {\tiny$$};
        \node at (7.1,-0.6) {\tiny$g$};
        \node at (6.2,-0.2) {\tiny$\vdots$};
        \node at (7.8,-0.2) {\tiny$\vdots$};
        \node at (4.5,0) {$G:$};
        \end{tikzpicture}
\caption{Local replacement in $\Sigma'$ of the special edge $g'$ in $G'$  and the resulting local picture in $O$ with new edge $g$  connecting the new marked point $X=Y$ with a new  orbifold point $\omega$. }
\label{fig:local-replacement}
\end{figure}
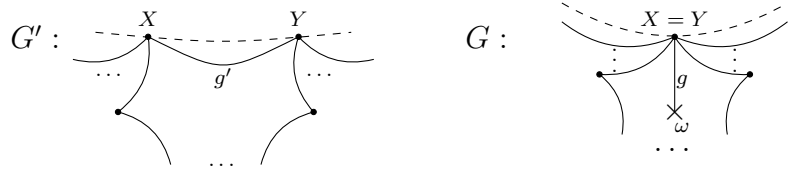
We do this for every edge corresponding to the start of a special loop and we obtain in this way a dissection of an orbifold which we will denote by  \sloppy $(O, \G, \G^*, M, P^*, \Omega)$, where $O$ is an orbifold with $\Omega$ its set of orbifold points (which are all of order two), $\G$ is a dissection of $O$ into  polygons (including the edges to orbifold points), $\G^*$ is the graph dual to $\G$, $M$ is the set of boundary vertices of $\G$, and $P^*$ are the vertices of $\G^*$ in the interior of $O$. 
By construction we have a bijection between $\Omega$ and the set of special loops $S$ and we call a polygon in $G$ containing at least one orbifold point a \emph{generalised polygon}. 

We then show   the following result, which has also been shown in \cite{AB}. 

\begin{proposition} With the notation above, 
there is a bijection between skew-gentle algebras up to isomorphism and homeomorphism classes of dissections of orbifolds into polygons and generalised polygons such that each (generalised) polygon has either exactly one boundary segment with a marked point in $M$ or is a (generalised) polygon with only internal edges with exactly one marked point  in $P^*$ in its interior. The bijection is given by the map which associates to a skew-gentle algebra $KQ/I$ the orbifold dissection $(O, \G, \G^*, M, P^*, \Omega)$.
\end{proposition}

We show that  a skew-gentle algebra $A$ is  Koszul and that its Koszul dual is skew-gentle with associated orbifold dissection induced by the dual graph of the orbifold dissection  of $A$.  More precisely, we show the following: 
 
 \begin{proposition}
 Let $A$ be a finite dimensional skew-gentle algebra. Then $A$ is Koszul and the Koszul dual $A^!$ is skew-gentle. Furthermore,  if $A$  has orbifold dissection $(O, \G, \G^*, M, P^*, \Omega)$ then $A^!$  has  orbifold dissection $(\tilde{O},  G^*, G, M^*,  P, \Omega)$ where  $\tilde{O}$ is the orbifold associated to $A!$,  $M^*$ are the boundary vertices of $G^*$ and $P$ corresponds to the set of non-boundary vertices of $\G$. 
 \end{proposition}

We note that similarly to the case of gentle algebras, a skew-gentle algebra $A$ is finite dimensional if and only if all the vertices of $\G$ are on the boundary of $O$, that is $P = \emptyset$. Furthermore, $A$ is of finite global dimension if and only if $P^* = \emptyset$.

Let $O$ be an orbifold as above. We  recall the notion of $O$-free homotopy from~\cite{CG}.

\begin{figure}[ht!]
     \centering
    \hspace*{0cm}    \begin{tikzpicture}[thick,scale=0.65, every node/.style={scale=0.65}]
       \draw (0,0) circle (1.5cm);
       \node at (0,0) {$\times$};
        \draw (-1.5,0) to[out=55, in=95] (0.75,0);
        \draw[<-] (-1.5,0) to[out=-55, in=-120] (0.35,-0.2);
        \draw(-0.35, -0.2) to[out=-65, in=-95] (0.75,0);
        \draw (-0.35,-0.2) .. controls (-0.4, -0.025) and (-0.175,0.15) .. (0,0.15)
               .. controls (0.175, 0.15) and (0.4,0.025) .. (0.35,-0.2);
        \node at (2,0) {$\sim$};
       \draw (4,0) circle (1.5cm);
       \node at (4,0) {$\times$};
       \draw[<-] (2.5,0) to[out=55, in=95] (3.5,0);
      \draw (2.5, 0) to[out=-55, in=-95] (3.5,0);
     \node at (6,0) {$\sim$};
       \draw (8,0) circle (1.5cm);
       \node at (8,0) {$\times$};
        \draw[<-](6.5,0) to[out=55, in=95] (8.75,0);
        \draw (6.5,0) to[out=-55, in=-120] (8.35,-0.2);
        \draw(7.65, -0.2) to[out=-65, in=-95] (8.75,0);
        \draw (7.65,-0.2) .. controls (7.6, -0.025) and (7.825,0.15) .. (8,0.15)
               .. controls (8.175, 0.15) and (8.4,0.025) .. (8.35,-0.2);
       \end{tikzpicture}
         \caption{Skein relations.}
    \label{fig:skein_relation2}
\end{figure}
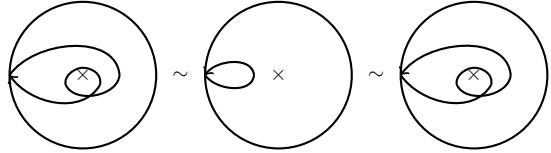

\begin{defn}
Two  oriented closed curves $\gamma$ and $\gamma'$ in an orbifold $O$ with orbifold points of order 2 are $O$-homotopic if they are related by a finite number of moves given by either a homotopy in the complement of the orbifold points or  a skein relation as in Figure~\ref{fig:skein_relation2} taking place in a disk $D$ containing exactly one orbifold point $\omega$. A segment of a curve with no self-intersection in $D$ and passing through $\omega$ is $O$-homotopic relative to its endpoints to a segment going around $\omega$ in either direction exactly once. 
\end{defn}

Let $\mathbb{H}$ denote the upper half plane and let $\Gamma <  PSL(2,\mathbb{R})$ be a Fuchsian group such that $O$ is a suborbifold of  $ \mathbb{H}/\Gamma $ with geodesic boundary. It is shown in \cite{CG} that there is a natural bijection between the set of conjugacy classes of $\Gamma$ and
the set of $O$-free homotopy classes of closed oriented curves in $O$. 

We are now in a position to state the main result of \cite{LSV}. 
\begin{theorem}
Let $A$ be a finite-dimensional skew-gentle algebra with orbifold dissection $(O, \G, \G^*, M, P^*, \Omega)$. Then the homotopy strings and bands of $A$, giving rise to the indecomposable objects in $D^b(A)$, are in bijection with graded curves $(\gamma, f)$ where 
\begin{enumerate}
\item $\gamma$ is an $O$-homotopy class of curves in  $O$  connecting marked points in $M \cup P^*$  or $\gamma$ is an $O$-homotopy class of certain closed curves in $O$. 
\item Given a curve $\gamma$ as in (1), set $f : \gamma \cap G^*  \to \mathbb{Z}$ to be the map such that, denoting the segment  connecting  two consecutive (in the direction of $\gamma$) intersection points $ x_i, x_{i+1}$ of $\gamma$ with $G^*$ by $\gamma_i$,  we have 
$$f( x_{ i+1}) = \left\{ 
\begin{array}{ll}
f(x_i)+1 & \mbox{if } \exists m \in M  \mbox{ to the left of $\gamma_i$ in $R_i$}  \\
f(x_i) -1 &  \mbox{if } \exists m \in M  \mbox{ to the right of $\gamma_i$ in  $R_i$} \\
\end{array} \right.
$$
where $R_i$ is the (generalised) polygon of $G^*$ containing $\gamma_i$.
\item  The open curves correspond to homotopy strings and the closed curves corresponding to homotopy bands are exactly those that have combinatorial winding number induced by $f$  equal to zero. 
\end{enumerate}
\end{theorem}

\begin{remark} {\rm 

(1) We note that the segment $\gamma_i$ will lie in exactly one (generalised) polygon of $G^*$ which by construction has exactly one boundary segment with exactly one marked point in $M$. 

(2) 
In ongoing work we show the connection between 'well-graded' intersections of two graded curves and maps between the corresponding objects in $D^b(A)$. }
\end{remark}

\end{document}